\documentclass{amsart}
\numberwithin{equation}{section}

\usepackage{amsmath,amssymb,mathtools}
\usepackage{booktabs,longtable,array}
\usepackage{enumitem}
\usepackage{microtype}

\newtheorem{theorem}{Theorem}[section]
\newtheorem{lemma}[theorem]{Lemma}

\begin{document}

\author{Ravshan Ashurov}
\address{Ashurov R: Institute of Mathematics of Uzbekistan,
Tashkent, Uzbekistan}
\email{ashurovr@gmail.com}

\author{Ilyoskhuja Sulaymonov}
\address{Sulaymonov I: National University of Uzbekistan,  Tashkent, Uzbekistan}
\email{ilyosxojasulaymonov@gmail.com}

\small

\title[Classical solvability of the Cauchy problem ...] {Classical solvability of the Cauchy problem for multi-order time-fractional systems with full matrices of differential operators}

\begin{abstract}

This article investigates the Cauchy problem in the whole space for a system of linear time-fractional partial differential equations with an arbitrary fixed number of components and a full matrix of constant-coefficient spatial differential operators. Each component is governed by its own Caputo derivative, and no rational relation between the fractional orders is required. We assume that the diagonal spatial symbols are homogeneous, positive and elliptic, while the off-diagonal entries have strictly lower order than the diagonal operator in the corresponding column. Using spatial Fourier and temporal Laplace transforms, we construct solution operators through convergent series of time convolutions involving scalar Mittag--Leffler kernels. These series account for all coupling paths, including repeated interactions between components. A fractional Duhamel representation incorporates nonhomogeneous source terms. We establish uniform convergence on bounded frequency regions and weighted multiplier estimates by separating low and high spatial frequencies. Under suitable Sobolev regularity assumptions on the initial data and continuously differentiable source terms, we prove existence and uniqueness of a classical solution, together with continuity and regularity properties in the appropriate Sobolev spaces. The general result includes two-component systems and is illustrated by fully coupled examples.

\end{abstract}

\maketitle

\section*{Introduction and problem statement}\label{sec:intro}

Systems with different fractional time orders describe interacting
quantities with different temporal responses. They also pose an analytic
difficulty: a change of dependent variables that diagonalizes a spatial
matrix need not commute with the diagonal matrix of time-fractional
derivatives. Consequently, scalar representations and estimates must
be combined without assuming a common fractional order.

For rational orders one can recover a common time derivative by
enlarging the system. Write $\beta_j=p_j/q_j$ in lowest terms and let $M$ be the
least common multiple of the denominators. The standard reduction gives
a system of order $1/M$ with $\sum_{j=1}^m M\beta_j$ components.
Umarov \cite{UmarovNew} illustrates the cost of this reduction with the
four orders $1/2$, $2/3$, $1/5$, and $6/7$. Here $M=210$, so the four
equations give rise to $105+140+42+180=467$ equations of order $1/210$.
For orders with an irrational ratio, even this finite reduction is
unavailable. For instance, if $1/2$ and $1/\sqrt{2}$ were integer
multiples of the same positive order, their ratio would be rational.
This raises the question of how to represent the solution directly in
terms of the original orders.

 In \cite{Umarovappear}, Umarov obtained matrix Mittag--Leffler
representations and uniqueness results for time-fractional
pseudo-differential systems, including particular cases with irrational
orders. He subsequently treated arbitrary vectors of orders for abstract
systems whose matrix entries are functions of a closed operator
\cite{UmarovNew}. The latter work constructs a matrix Mittag--Leffler
function indexed by the vector of orders and proves existence and
uniqueness in the associated topological vector spaces. For a partial
differential system on $\mathbb R^n$, a further question is which Sobolev
conditions on the data ensure that the represented solution satisfies
each equation pointwise.

Ashurov and Mukhiddinova \cite{AO} addressed classical solvability for
systems with a common Caputo order and a Hermitian spatial symbol.
They treated both a prescribed source and a source depending on the
solution. In the latter case they imposed a Lipschitz condition in the
unknown together with an additional hypothesis for global existence.
For independent Caputo orders, Ashurov and Sulaymonov \cite{AS} obtained
conditions for classical solvability when the spatial operator matrix
is triangular and its diagonal entries are elliptic. Triangularity allows the equations to be solved successively,
and the fractional Duhamel principle gives the contribution of the
prescribed source. Sulaymonov \cite{SulaymonovRL} studied the corresponding
Riemann--Liouville problem, for which the initial data are prescribed
through fractional integrals when $0<\beta_j<1$.

We return to the Caputo problem and ask whether classical solvability
persists when the spatial operator matrix is full. The difficulty is
already visible in a two-component system. If each equation contains
both unknowns, substituting the second component into the first brings
the first component back into its own equation. Successive substitution
then produces an infinite series of time convolutions. To obtain a
classical solution on $\mathbb R^n$, we shall estimate this series uniformly
in the spatial frequency.

Let $T>0$, $m\ge2$ and $\mathcal B=(\beta_1,\dots,\beta_m)$, where
$0<\beta_j\le1$. We consider
\begin{align}
 D_t^{\mathcal B}U(t,x)+\mathbb A(D)U(t,x)&=H(t,x),
 &&0<t\le T,\quad x\in\mathbb R^n,\label{eq:cauchy}\\
 U(0,x)&=\Phi(x),&&x\in\mathbb R^n,\label{eq:initial}
\end{align}
where
\[
 U=(u_1,\dots,u_m)^\mathsf T,\quad
 H=(h_1,\dots,h_m)^\mathsf T,\quad
 \Phi=(\varphi_1,\dots,\varphi_m)^\mathsf T,
D_t^{\mathcal B}U
 =(D_t^{\beta_1}u_1,\dots,D_t^{\beta_m}u_m)^\mathsf T.\]
Let $\mathbb A(D)=(A_{i,j}(D))_{i,j=1}^m$ be a matrix of
differential operators with constant coefficients, $ A_{i,j}(D)=\sum_{|\alpha|\leq \ell_{i,j}}a_{\alpha}^{(i,j)}D^\alpha,$
where $\ell_{i,j}$ denotes the order of $A_{i,j}(D)$. Thus
\[
\mathbb{A}(D)= \begin{bmatrix}
	    A_{1,1}(D)& A_{1,2}(D) & \dots & A_{1,m}(D) \\
	    A_{2,1}(D)& A_{2,2}(D) &\dots & A_{2,m}(D)\\
	    \dots  & \dots  & \dots & \dots  \\
	    A_{m, 1}(D) & A_{m,2}(D) & \dots & A_{m,m} (D)
	    \end{bmatrix},\] 
with the matrix-symbol
\[
\mathcal{A}(\xi) = \begin{bmatrix}
	    A_{1,1}(\xi)& A_{1,2}(\xi) & \dots & A_{1,m}(\xi) \\
	    A_{2,1}(\xi)& A_{2,2}(\xi) &\dots & A_{2,m}(\xi)\\
	    \dots  & \dots  & \dots & \dots  \\
	    A_{m, 1}(\xi) & A_{m,2}(\xi) & \dots & A_{m,m} (\xi)
	    \end{bmatrix},
\]
Here $\alpha=(\alpha_1,\ldots,\alpha_n)$ is a multi-index,
$D=(D_1,\ldots,D_n)$ with
$D_j=\frac{1}{\mathrm i}\frac{\partial}{\partial x_j}$, and
$ A_{i,j}(\xi)=\sum_{|\alpha|\leq \ell_{i,j}}
 a_{\alpha}^{(i,j)}\xi^\alpha$
is the symbol of $A_{i,j}(D)$.

Throughout the article we assume that the diagonal symbols $A_{jj}(\xi)$
are real, homogeneous, elliptic, and positive for $\xi\ne0$:
\begin{equation}\label{eq:positive}
 c_j|\xi|^{\ell_{jj}}\le A_{jj}(\xi)
 \le C_j|\xi|^{\ell_{jj}},\quad \xi\ne0,
 \qquad j=1,\dots,m,
\end{equation}
where each $\ell_{jj}$ is a positive even integer. Moreover,
\begin{equation}\label{eq:column}
 \ell_{ij}<\ell_{jj},\qquad i\ne j,\quad i,j=1,\dots,m.
\end{equation}
Thus the diagonal operator has the highest order in each column.

For $s\in\mathbb R$, let $L_2^s(\mathbb R^n)$ be the Bessel-potential Sobolev space
with norm
\[
 \|f\|_{L_2^s}^2=(2\pi)^{-n}
 \int_{\mathbb R^n}\langle \xi\rangle^{2s}|\widehat f(\xi)|^2\,\mathrm d\xi,
 \qquad \langle \xi\rangle=(1+|\xi|^2)^{1/2}.
\]
Write $\mathbf L_2^\tau=(L_2^\tau)^m$,
$\mathbf C(E)=C(E;\mathbb C^m)$ for a space-time set $E$, and introduce
\[
 \mathbf X_\tau=\prod_{j=1}^m L_2^{\tau+\ell_{jj}}(\mathbb R^n),
 \quad
 \|U\|_{\mathbf X_\tau}^2
 =\sum_{j=1}^m\|u_j\|_{L_2^{\tau+\ell_{jj}}}^2.
\]

For $0<\beta<1$ we use the Caputo derivative
\[
 D_t^\beta f(t)=\frac{\,\mathrm d}{\,\mathrm d t}
 I_t^{1-\beta}(f-f(0))(t),
 \qquad
 I_t^\gamma f(t)=\frac1{\Gamma(\gamma)}
 \int_0^t(t-s)^{\gamma-1}f(s)\,\mathrm d s.
\]
For $\beta=1$, we set $D_t^1f=f'$.

\noindent\textbf{Cauchy problem and solution class.}
\textit{Find functions $u_j(t,x)\in L_2^{\ell_{jj}}(\mathbb R^n)$,
$t\in(0,T]$, $j=1,\dots,m$, satisfying \eqref{eq:cauchy}--\eqref{eq:initial}
and}
\begin{equation}\label{eq:classical}
 U\in\mathbf C([0,T]\times\mathbb R^n),\qquad
 D_t^{\mathcal B}U,\ \mathbb A(D)U
 \in\mathbf C((0,T]\times\mathbb R^n).
\end{equation}
\textit{Such a solution will be called a classical solution in the $L_2$ class.}
Here $H$ and $\Phi$ are prescribed continuous vector functions.
The functional restriction at spatial infinity is expressed through the
$L_2$ regularity. 

\section{Preliminaries}\label{sec:prelim}

The Fourier transform and its inverse are defined by
\[
 \mathrm F[f](\xi)=\int_{\mathbb R^n}e^{-\mathrm ix\cdot\xi}f(x)\,\mathrm d x,
 \qquad
 \mathrm F^{-1}[g](x)=(2\pi)^{-n}\int_{\mathbb R^n}e^{\mathrm ix\cdot\xi}g(\xi)\,\mathrm d\xi.
\]
In particular, $\mathrm F[A_{ij}(D)f]=A_{ij}(\xi)\widehat f(\xi)$.
We use the time convolution
\[
 (f*g)(t)=\int_0^tf(t-s)g(s)\,\mathrm d s,
\]
and $f^{*r}$ denotes its $r$th convolution power. For matrices,
$(F*G)_{ij}=\sum_kF_{ik}*G_{kj}$; the order of matrix factors is retained.
The zeroth convolution power is the identity convolution operator.

The two-parameter Mittag--Leffler function is
\[
 E_{\beta,\mu}(z)=\sum_{k=0}^\infty\frac{z^k}{\Gamma(\beta k+\mu)},
 \qquad E_\beta(z)=E_{\beta,1}(z).
\]
At every fixed frequency $\xi$, put $A_{jj}(\xi)$ and define
\[
 S_j(t,\xi)=E_{\beta_j}(-A_{jj}(\xi)t^{\beta_j}),\quad
 K_j(t,\xi)=t^{\beta_j-1}
 E_{\beta_j,\beta_j}(-A_{jj}(\xi)t^{\beta_j}),
 \qquad j=1,\dots,m.
\]
For $A_{j,j}\ge0$ and $0<\beta_j\le1$, standard scalar estimates
(see \cite{M,Gor}) give
\begin{equation}\label{eq:scalar-est}
 0\le S_j(t,\xi)\le1,\quad
 0\le K_j(t,\xi)\le\frac{t^{\beta_j-1}}{\Gamma(\beta_j)}.
\end{equation}
For $\xi\ne0$, we have
\[
 \int_0^T K_j(t,\xi)\,\mathrm d t
 =\frac{1-S_j(T,\xi)}{A_{jj}(\xi)}
 \le\frac1{A_{jj}(\xi)}.
\]
The Laplace transforms are
\begin{equation}\label{eq:scalar-L}
 \mathrm L[S_j](s,\xi)=\frac{s^{\beta_j-1}}{s^{\beta_j}+A_{jj}(\xi)},
 \qquad
 \mathrm L[K_j](s,\xi)=\frac1{s^{\beta_j}+A_{jj}(\xi)}.
\end{equation}

For $\alpha>0$, set $g_\alpha(t)=t^{\alpha-1}/\Gamma(\alpha)$.
The Beta integral yields
\begin{equation}\label{eq:beta-conv}
 g_\alpha*g_\gamma=g_{\alpha+\gamma},\qquad
 (g_\alpha*1)(t)=\frac{t^\alpha}{\Gamma(\alpha+1)}.
\end{equation}
Finally, since $\tau>n/2$, the Sobolev embedding and Fourier--H\"older
inequality give
\begin{equation}\label{eq:embedding}
 L_2^\tau(\mathbb R^n)\hookrightarrow C_b(\mathbb R^n),\qquad
 \|f\|_{C_b(\mathbb R^n)}\le C_\tau\|f\|_{L_2^\tau}.
\end{equation}

\section{The two-component case}\label{sec:two}
We first study the case $m=2$ in detail. The calculations below will also
indicate the estimates that are needed for the general case. We consider
\begin{align}
 D_t^{\beta_1}u_1+A_{1,1}(D)u_1+A_{1,2}(D)u_2&=h_1,\label{tw:eq1}\\
 D_t^{\beta_2}u_2+A_{2,1}(D)u_1+A_{2,2}(D)u_2&=h_2,\label{tw:eq2}\\
 u_1(0,x)=\varphi_1(x),\qquad u_2(0,x)&=\varphi_2(x).\label{tw:initial}
\end{align}
The diagonal symbols $A_{jj}(\xi)$ satisfy \eqref{eq:positive},
and the column-order conditions become
\begin{equation}\label{tw:orders}
 \ell_{11}>\ell_{21},\qquad \ell_{22}>\ell_{12}.
\end{equation}

\subsection{Formal solution}\label{sec:two-formal}
We first put $h_1=h_2=0$. Taking the Fourier transform of
\eqref{tw:eq1}--\eqref{tw:initial} with respect to $x$ gives
\[
\left\{
\begin{aligned}
    & D_t^{\mathcal{B}} \mathrm F[U](t,\xi) +\mathbb{A}(\xi) \mathrm F[U](t,\xi)=0, \quad t>0, \ \xi \in \mathbb{R}^n,\\
&\mathrm F[U](0,\xi) = \mathrm F[\varPhi](\xi), \quad \xi \in \mathbb{R}^n.
\end{aligned}\right.
\]
Then apply the Laplace transform to get
\begin{equation}\label{tw:algebraic}
\begin{bmatrix}
s^{\beta_1} + A_{1,1}(\xi)& A_{1,2}(\xi)  \\
A_{2,1}(\xi) & s^{\beta_2}+A_{2,2}(\xi)
\end{bmatrix}
\begin{bmatrix}
	\mathrm L[\mathrm F[u_1]](s,\xi)\\
	\mathrm L[\mathrm F[u_2]](s,\xi)
\end{bmatrix}
=
\begin{bmatrix}
	s^{\beta_1-1}\mathrm F[\varphi_1](\xi)\\
	s^{\beta_2-1}\mathrm F[\varphi_2](\xi)
\end{bmatrix}.
\end{equation}
Let
\[
 \Delta(s,\xi)=(s^{\beta_1}+A_{1,1}(\xi))(s^{\beta_2}+A_{2,2}(\xi))-A_{1,2}(\xi)A_{2,1}(\xi), \quad P(\xi)=A_{1,2}(\xi)A_{2,1}(\xi).
\]
By solving \eqref{tw:algebraic}, we obtain
\[
\begin{aligned}
\mathrm L[\mathrm F[u_1]](s,\xi)&=\frac{(s^{\beta_2}+A_{2,2})s^{\beta_1-1}\mathrm F[\varphi_1](\xi)
 -A_{1,2}s^{\beta_2-1}\mathrm F[\varphi_2](\xi)}{\Delta(s,\xi)},\\
 \mathrm L[\mathrm F[u_2]](s,\xi)&=\frac{-A_{2,1}s^{\beta_1-1}\mathrm F[\varphi_1](\xi)
 +(s^{\beta_1}+A_{1,1})s^{\beta_2-1}\mathrm F[\varphi_2](\xi)}{\Delta(s,\xi)}.
\end{aligned}
\]
For sufficiently large real $s$ at each fixed $\xi$,
\[
\frac1{\Delta(s,\xi)}=\sum_{k=0}^\infty\frac{P^k}
 {(s^{\beta_1}+A_{1,1})^{k+1}(s^{\beta_2}+A_{2,2})^{k+1}}.
\]
Formula \eqref{eq:scalar-L} and the convolution equality yield
\begin{align}
 S_{11}(t,\xi)&=\sum_{k=0}^\infty P^k
 (K_1^{*k}*K_2^{*k}*S_1)(t,\xi),\label{tw:S11}\\
 S_{12}(t,\xi)&=-A_{12}\sum_{k=0}^\infty P^k
 (K_1^{*(k+1)}*K_2^{*k}*S_2)(t,\xi),\label{tw:S12}\\
 S_{21}(t,\xi)&=-A_{21}\sum_{k=0}^\infty P^k
 (K_1^{*k}*K_2^{*(k+1)}*S_1)(t,\xi),\label{tw:S21}\\
 S_{22}(t,\xi)&=\sum_{k=0}^\infty P^k
 (K_1^{*k}*K_2^{*k}*S_2)(t,\xi).\label{tw:S22}
\end{align}
Therefore the formal homogeneous solution takes the form
\[
 U(t,x)=S(t,D)\Phi(x),\qquad
 u_i(t,x)=\frac1{(2\pi)^n}\sum_{j=1}^2\int_{\mathbb R^n}
 e^{\mathrm ix\cdot\xi}S_{ij}(t,\xi)\widehat\varphi_j(\xi)\,\mathrm d\xi.
\]

For $H\ne0$ we apply the fractional Duhamel principle
\cite{Umarov1,Umarov2}. We denote the source kernel by $S'$; the prime
\emph{does not mean ordinary differentiation in $t$}.
The kernels are
\begin{align}
 S'_{11}(t,\xi)&=\sum_{k=0}^\infty P^k
 (K_1^{*(k+1)}*K_2^{*k})(t,\xi),\label{tw:Sp11}\\
 S'_{12}(t,\xi)&=-A_{12}\sum_{k=0}^\infty P^k
 (K_1^{*(k+1)}*K_2^{*(k+1)})(t,\xi),\label{tw:Sp12}\\
 S'_{21}(t,\xi)&=-A_{21}\sum_{k=0}^\infty P^k
 (K_1^{*(k+1)}*K_2^{*(k+1)})(t,\xi),\label{tw:Sp21}\\
 S'_{22}(t,\xi)&=\sum_{k=0}^\infty P^k
 (K_1^{*k}*K_2^{*(k+1)})(t,\xi).\label{tw:Sp22}
\end{align}
Thus the formal solution to the nonhomogeneous problem is
\begin{equation}\label{tw:solution}
 U(t)=S(t,D)\Phi+\int_0^tS'(t-s,D)H(s)\,\mathrm d s.
\end{equation}

\subsection{Convergence and estimates of the four factors}
\label{sec:two-bounds}

We derive the estimates for the solution factors
$S_{ij}$ and the source kernels $S'_{ij}$.

\begin{lemma}\label{lem:two-bounds}
Assume \eqref{eq:positive} and \eqref{tw:orders}.
Then the series \eqref{tw:S11}--\eqref{tw:S22}
converge absolutely and uniformly on every set
$[0,T]\times\{|\xi|\le R\}$.

Moreover, the source kernels $S'_{ij}$ are
integrable in time, and for
$t\in[0,T]$, $\xi\in\mathbb R^n$, and $i,j=1,2$,
the following estimates hold:
\begin{align}
 \langle \xi\rangle^{\ell_{ii}}|S_{ij}(t,\xi)|
 &\le C_T\langle \xi\rangle^{\ell_{jj}},
 \label{tw:S-weighted}\\
 \langle \xi\rangle^{\ell_{ii}}
 \int_0^T|S'_{ij}(s,\xi)|\,\mathrm d s
 &\le C_T.
 \label{tw:Sp-weighted}
\end{align}
\end{lemma}

\begin{proof}
\textit{1. The bounded-frequency region $|\xi|\le R$.}
Put $\sigma=\beta_1+\beta_2$ and
\[
a_R=\max_{|\xi|\le R}|A_{12}(\xi)|,\quad
b_R=\max_{|\xi|\le R}|A_{21}(\xi)|,\quad p_R=a_Rb_R.
\]
By \eqref{eq:scalar-est} and \eqref{eq:beta-conv},
$0\le S_j\le1$, $0\le K_j\le g_{\beta_j}$, and
\[
(K_1^{*k}*K_2^{*l}*S_j)(t,\xi)
\le\frac{t^{k\beta_1+l\beta_2}}
{\Gamma(k\beta_1+l\beta_2+1)}.
\]
Substitution into \eqref{tw:S11}--\eqref{tw:S22} yields
\[
\begin{aligned}
|S_{11}|,|S_{22}|&\le E_{\sigma,1}(p_Rt^\sigma),\\
|S_{12}|&\le a_Rt^{\beta_1}
E_{\sigma,\beta_1+1}(p_Rt^\sigma),\\
|S_{21}|&\le b_Rt^{\beta_2}
E_{\sigma,\beta_2+1}(p_Rt^\sigma).
\end{aligned}
\]
Since these majorants are finite on $[0,T]$, the four
series converge absolutely and uniformly for $|\xi|\le R$.
Similarly, integrating the convolution estimates for
\eqref{tw:Sp11}--\eqref{tw:Sp22} gives
\[
\begin{aligned}
\int_0^t|S'_{11}(s,\xi)|\,\mathrm d s
&\le t^{\beta_1}E_{\sigma,\beta_1+1}(p_Rt^\sigma),\\
\int_0^t|S'_{22}(s,\xi)|\,\mathrm d s
&\le t^{\beta_2}E_{\sigma,\beta_2+1}(p_Rt^\sigma),\\
\int_0^t|S'_{12}(s,\xi)|\,\mathrm d s
&\le a_Rt^\sigma E_{\sigma,\sigma+1}(p_Rt^\sigma),\\
\int_0^t|S'_{21}(s,\xi)|\,\mathrm d s
&\le b_Rt^\sigma E_{\sigma,\sigma+1}(p_Rt^\sigma).
\end{aligned}
\]
Thus the source series converge absolutely in $L_1(0,T)$,
uniformly for $|\xi|\le R$.

\medskip
\textit{2. The large-frequency region $|\xi|\ge R$.}
By ellipticity and \eqref{tw:orders},
\[
q(\xi):=\frac{|A_{12}(\xi)A_{21}(\xi)|}
{A_{1,1}(\xi)A_{2,2}(\xi)}
\le C|\xi|^{\ell_{12}+\ell_{21}-\ell_{11}-\ell_{22}}
\longrightarrow0 \quad (|\xi|\to\infty).
\]
Choose $R>1$ so large that $q(\xi)\le1/2$ whenever
$|\xi|\ge R$. Using $0\le S_j\le1$ and
$\int_0^T K_j(s,\xi)\,\mathrm d s\le A_{j,j}(\xi)^{-1}$,
the series in \eqref{tw:S11}--\eqref{tw:S22} are bounded by
\[
|S_{11}|,|S_{22}|\le \sum_{k=0}^\infty q(\xi)^k\le2,\quad
|S_{12}|\le\frac{2|A_{12}(\xi)|}{A_{1,1}(\xi)},\quad
|S_{21}|\le\frac{2|A_{21}(\xi)|}{A_{2,2}(\xi)}.
\]
The same argument applied to the source kernels gives
\[
\begin{aligned}
\int_0^T|S'_{11}(s,\xi)|\,\mathrm d s
\le\frac2{A_{1,1}(\xi)},&\quad
\int_0^T|S'_{22}(s,\xi)|\,\mathrm d s
\le\frac2{A_{2,2}(\xi)},\\
\int_0^T|S'_{12}(s,\xi)|\,\mathrm d s
\le\frac{2|A_{12}(\xi)|}{A_{1,1}(\xi)A_{2,2}(\xi)},&\quad
\int_0^T|S'_{21}(s,\xi)|\,\mathrm d s
\le\frac{2|A_{21}(\xi)|}{A_{1,1}(\xi)A_{2,2}(\xi)}.
\end{aligned}
\]
For $|\xi|\ge R$, ellipticity implies
$A_{jj}(\xi)\asymp|\xi|^{\ell_{jj}}$, whereas
$|A_{ij}(\xi)|\le C|\xi|^{\ell_{ij}}$.
Consequently, the column-order inequalities
$\ell_{12}<\ell_{22}$ and $\ell_{21}<\ell_{11}$
yield \eqref{tw:S-weighted} and \eqref{tw:Sp-weighted}.
The bounds on $|\xi|\le R$ obtained in part 1 give
these estimates on the remaining frequency region as well.
\end{proof}

\subsection{Main theorem and its proof}\label{sec:two-main}

For the two-component system define
\[
\mathbf X_\tau^{(2)}
=L_2^{\tau+\ell_{11}}(\mathbb R^n)
\times L_2^{\tau+\ell_{22}}(\mathbb R^n).
\]

\begin{theorem}\label{th:two}
Let $\tau>n/2$, $0<\beta_1,\beta_2\le1$,
and assume \eqref{eq:positive} and
\eqref{tw:orders}. If
\[
\varphi_j\in L_2^{\tau+\ell_{jj}}(\mathbb R^n),
\qquad
h_j\in C^1([0,T];L_2^\tau(\mathbb R^n)),
\quad j=1,2,
\]
then \eqref{tw:eq1}--\eqref{tw:initial}
has a unique classical solution in
$C([0,T];\mathbf X_\tau^{(2)})$.
This solution is given by \eqref{tw:solution}.
\end{theorem}

\begin{proof}
\textit{Existence and continuity.}

Using \eqref{tw:S-weighted}, we obtain
\begin{equation}\label{tw:Sestimate}
\|S(t,D)\Phi\|_{\mathbf X_\tau^{(2)}}
\le C_T\|\Phi\|_{\mathbf X_\tau^{(2)}}.
\end{equation}

We now estimate the source term
\[
V(t)=\int_0^t S'(t-s,D)H(s)\,\mathrm d s.
\]

Since $H\in C^1([0,T];(L_2^\tau)^2)$,
we have
\[
H(s)=H(0)+\int_0^s H'(r)\,\mathrm d r.
\]
Consequently, for almost every $\xi$,
\begin{equation}\label{tw:H-Fourier}
|\mathrm F[h_j](s,\xi)|
\le |\mathrm F[h_j](0,\xi)|
+\int_0^T|\mathrm F[h_j'](r,\xi)|\,\mathrm d r.
\end{equation}

Taking the Fourier transform of the
source term gives
\[
\mathrm F[V_i](t,\xi)
=\sum_{j=1}^2
\int_0^t S'_{ij}(t-s,\xi)
\mathrm F[h_j](s,\xi)\,\mathrm d s.
\]

Using \eqref{tw:H-Fourier} and
\eqref{tw:Sp-weighted}, we obtain
\[
\begin{aligned}
&\langle \xi\rangle^{\ell_{ii}}
|\mathrm F[V_i](t,\xi)|\le
\sum_{j=1}^2
\langle \xi\rangle^{\ell_{ii}}
\int_0^t|S'_{ij}(t-s,\xi)|
|\mathrm F[h_j](s,\xi)|\,\mathrm d s\\
&\le C_T\sum_{j=1}^2
\left(
|\mathrm F[h_j](0,\xi)|
+\int_0^T|\mathrm F[h_j'](r,\xi)|\,\mathrm d r
\right).
\end{aligned}
\]

Multiplying by $\langle \xi\rangle^{\tau}$,
integrating over $\mathbb R^n$, and applying
Minkowski's inequality, we find
\begin{align}
&\left\|
\int_0^t S'(t-s,D)H(s)\,\mathrm d s
\right\|_{\mathbf X_\tau^{(2)}}
\le C_T
\left(
\|H(0)\|_{(L_2^\tau)^2}
+\int_0^T
\|H'(s)\|_{(L_2^\tau)^2}\,\mathrm d s
\right)
\notag\\
&\quad\le C_T
\left(
\|H(0)\|_{(L_2^\tau)^2}
+T\|H'\|_{C([0,T];(L_2^\tau)^2)}
\right).
\label{tw:sourceestimate}
\end{align}

Combining \eqref{tw:Sestimate} and
\eqref{tw:sourceestimate}, we obtain
\[
\begin{aligned}
&\sup_{0\le t\le T}
\|U(t)\|_{\mathbf X_\tau^{(2)}}\le C_T\bigl(
\|\Phi\|_{\mathbf X_\tau^{(2)}}
+\|H(0)\|_{(L_2^\tau)^2}
+T\|H'\|_{C([0,T];(L_2^\tau)^2)}
\bigr).
\end{aligned}
\]

The preceding Fourier estimates
provide an integrable majorant
independent of $t$.
Hence 
\begin{equation}\label{tw:cont}
U\in C([0,T];\mathbf X_\tau^{(2)}).
\end{equation}

\textit{Continuity of the spatial expressions.}

By \eqref{tw:orders}, for $i,j=1,2$,
we have $\ell_{ij}\le\ell_{jj}$.
Hence
\[
\begin{aligned}
\|A_{ij}(D)u_j(t)\|_{L_2^\tau}
&\le C\|u_j(t)\|_{L_2^{\tau+\ell_{ij}}}.
\end{aligned}
\]

Together with \eqref{tw:cont},
this proves
\[
\mathbb A(D)U
\in C([0,T];(L_2^\tau)^2).
\]

Since $\tau>n/2$, the Sobolev
embedding theorem gives
\[
U\in\mathbf C([0,T]\times\mathbb R^n),
\qquad
\mathbb A(D)U
\in\mathbf C((0,T]\times\mathbb R^n).
\]

Moreover,
$u_j(t,x)\in
L_2^{\ell_{jj}}(\mathbb R^n)$.

\medskip
\textit{Verification of the equations and the initial condition.}
The construction of $S$ and $S'$ gives, for almost every $\xi$,
\begin{equation}\label{tw:full-recurrence}
 \widehat u_i
 =S_i\widehat\varphi_i
 +K_i*\left(\widehat h_i-
 \sum_{\substack{j=1\\j\ne i}}^2 A_{ij}(\xi)\widehat u_j\right),
 \qquad i=1,2.
\end{equation}
Hence
\eqref{tw:full-recurrence} implies
\[
 D_t^{\mathcal B}\widehat U(t,\xi)+\mathbb A(\xi)\widehat U(t,\xi)
 =\widehat H(t,\xi)
\]
for $t>0$ and almost every $\xi$. Since both sides belong to
$(L_2^\tau)^2$, inversion of the Fourier transform yields
\eqref{tw:eq1}--\eqref{tw:eq2} in $(L_2^\tau)^2$.

The weighted
bounds of Lemma~\ref{lem:two-bounds} provide an $L_2$ majorant, and
therefore dominated convergence gives
\[
 U(t)\longrightarrow\Phi\quad\text{in }\mathbf X_\tau^{(2)}
 \qquad(t\rightarrow0).
\]
Thus \eqref{tw:initial} holds. Finally, from the equation,
\[
 D_t^{\mathcal B}U=H-\mathbb A(D)U
 \in C((0,T];(L_2^\tau)^2),
\]
and \eqref{eq:embedding} gives the corresponding spatial continuity.
Hence the function constructed above is a classical solution in the
sense of \eqref{eq:classical}.

\medskip
\textit{Uniqueness.}
Let $U^{(1)}$ and $U^{(2)}$ be two classical
solutions with the same initial data and
right-hand side. Then $W=U^{(1)}-U^{(2)}$
satisfies
\[
D_t^{\mathcal B}W+\mathbb A(D)W=0,
\qquad W(0)=0.
\]

Taking the Fourier transform with respect
to $x$, we obtain
\[
D_t^{\mathcal B}\mathrm F[W](t,\xi)
+\mathbb A(\xi)\mathrm F[W](t,\xi)=0,
\qquad \mathrm F[W](0,\xi)=0.
\]

For almost every fixed $\xi$, this is a
homogeneous linear system of fractional
ordinary differential equations with
constant coefficients and zero initial data.
By the standard uniqueness result for
linear multi-order fractional systems
(see \cite{UmarovNew}), we obtain
\[
\mathrm F[W](t,\xi)=0.
\]

Consequently, by continuity in time, $W=0$.
Hence $U^{(1)}=U^{(2)}$, and the
solution is unique.

\end{proof}

\section{The general case}\label{sec:general}
We now extend the preceding construction to any fixed $m\ge2$.
Throughout this section we retain the same scalar kernels $S_j$, $K_j$,
the solution factors $S_{ij}$ and the source kernels $S'_{ij}$ as in
the two-component case. In particular, the prime in $S'$ does not
denote differentiation with respect to $t$.

\subsection{Formal solution}\label{sec:formal}
We first consider the homogeneous problem $H=0$. Applying the Fourier
transform in $x$ and then the Laplace transform in $t$ gives
\[
 \bigl[\operatorname{diag}(s^{\beta_1},\ldots,s^{\beta_m})
       +\mathbb A(\xi)\bigr]
 \mathrm L[\mathrm F[U]](s,\xi)
 =\operatorname{diag}(s^{\beta_1-1},\ldots,s^{\beta_m-1})
 \mathrm F[\Phi](\xi).
\]
Equivalently, for every $i=1,\ldots,m$,
\[
 (s^{\beta_i}+A_{i,i}(\xi))\mathrm L[\mathrm F[u_i]](s,\xi)=s^{\beta_i-1}\mathrm F[\varphi_i](\xi)
 -\sum_{\substack{j=1\\j\ne i}}^{m}
A_{ij}(\xi)\mathrm L[\mathrm F[u_j]](s,\xi).
\]
Using the scalar transforms \eqref{eq:scalar-L}, we obtain

\begin{equation}\label{eq:general-recurrence}
\mathrm F[u_i](t,\xi)
={}S_i(t,\xi)\mathrm F[\varphi_i](\xi)-\sum_{\substack{j=1\\j\ne i}}^{m}
A_{ij}(\xi)
(K_i*\mathrm F[u_j])(t,\xi),
\quad i=1,\ldots,m.
\end{equation}
Here convolution acts only on the time variable. Successively
substituting \eqref{eq:general-recurrence} into its right-hand side,
exactly as in the two-component case, yields
\begin{equation}\label{eq:Sij}
\begin{aligned}
S_{ij}(t,\xi)
={}\delta_{ij}S_i(t,\xi)+\sum_{r=1}^{\infty}(-1)^r
\sum_{\substack{
i_0=i,\ i_r=j\\
i_{q-1}\ne i_q,\ 1\le q\le r}}
\left(
\prod_{q=1}^{r}
A_{i_{q-1},i_q}(\xi)
\right)
(K_{i_0}*\cdots*K_{i_{r-1}}*S_j)(t,\xi).
\end{aligned}
\end{equation}
The intermediate indices $i_1,\ldots,i_{r-1}$ range over
$\{1,\ldots,m\}$; $\delta_{ij}$ is the Kronecker symbol.
Thus $S(t,\xi)=(S_{ij}(t,\xi))_{i,j=1}^m$ and the formal
homogeneous solution is
\[
 U(t,x)=S(t,D)\Phi(x),\quad
 u_i(t,x)=\frac1{(2\pi)^n}\sum_{j=1}^m\int_{\mathbb R^n}
 e^{\mathrm ix\cdot\xi}S_{ij}(t,\xi)
 \mathrm F[\varphi_j](\xi)\,\mathrm d\xi.
\]

For $H\ne0$, applying the fractional Duhamel principle
\cite{Umarov1,Umarov2}, or repeating the same substitution with
the additional term $K_i*\mathrm F[h_i]$, gives the source kernels
\begin{equation}\label{eq:Spij}
 \begin{aligned}
 S'_{ij}(t,\xi)={}\delta_{ij}K_i(t,\xi)+\sum_{r=1}^{\infty}(-1)^r
 \sum_{\substack{i_0=i,\ i_r=j\\
                     i_{q-1}\ne i_q,\ 1\le q\le r}}
 \left(\prod_{q=1}^{r}A_{i_{q-1},i_q}(\xi)\right)
 (K_{i_0}*\cdots*K_{i_{r-1}}*K_j)(t,\xi).
 \end{aligned}
\end{equation}
Put $S'(t,\xi)=(S'_{ij}(t,\xi))_{i,j=1}^m$.
In particular, for sufficiently large real $s$ at each fixed $\xi$,
\[
 \mathrm L[S](s,\xi)
 =\bigl[\operatorname{diag}(s^{\beta_j})
     +\mathbb A(\xi)\bigr]^{-1}
   \operatorname{diag}(s^{\beta_j-1}),\quad
 \mathrm L[S'](s,\xi)
 =\bigl[\operatorname{diag}(s^{\beta_j})
     +\mathbb A(\xi)\bigr]^{-1}.
 \]
These formulas agree with the four factors and four source kernels
of Section~\ref{sec:two} when $m=2$. Hence the candidate solution is
\begin{equation}\label{eq:representation}
 U(t,x)=S(t,D)\Phi(x)
 +\int_0^tS'(t-s,D)H(s,x)\,\mathrm d s.
\end{equation}

\subsection{Convergence and estimates of the solution factors}
\label{sec:bounds}

\begin{lemma}\label{lem:multipliers}
Assume \eqref{eq:positive} and \eqref{eq:column}.
Then the series \eqref{eq:Sij} converge absolutely and uniformly on
every set $[0,T]\times\{|\xi|\le R\}$, and the source series
\eqref{eq:Spij} converge in $L_1(0,T)$ uniformly for $|\xi|\le R$.
Moreover, for $i,j=1,\ldots,m$, $0\le t\le T$ and $\xi\in\mathbb R^n$,
\begin{align}
 \langle \xi\rangle^{\ell_{ii}}|S_{ij}(t,\xi)|
 &\le C_T\langle \xi\rangle^{\ell_{jj}},\label{eq:S-bound}\\
 \langle \xi\rangle^{\ell_{ii}}
 \int_0^T|S'_{ij}(s,\xi)|\,\mathrm d s
 &\le C_T.\label{eq:Sp-bound}
\end{align}
\end{lemma}

\begin{proof}
\textit{1. The bounded-frequency region $|\xi|\le R$.}
Fix $R>0$ and set
\[
 \beta_*=\min_{1\le j\le m}\beta_j,
 \qquad
 \rho_R=\max_{\substack{|\xi|\le R\\i\ne j}}
 |A_{ij}(\xi)|.
\]
Since $0<\beta_*\le\beta_j\le1$, we can choose
\[
 c_{\mathcal B,T}=\max_{1\le j\le m}
 \frac{\Gamma(\beta_*)}{\Gamma(\beta_j)}
 T^{\beta_j-\beta_*}
\]
so that, for $0<s\le T$,
\[
 0\le K_j(s,\xi)\le
 c_{\mathcal B,T}g_{\beta_*}(s),
 \qquad 0\le S_j(s,\xi)\le1.
\]
By the Beta identity \eqref{eq:beta-conv}, for every chain
$i_0,\ldots,i_r$ in \eqref{eq:Sij},
\[
 (K_{i_0}*\cdots*K_{i_{r-1}}*S_j)(t,\xi)\le c_{\mathcal B,T}^{r}
 \frac{t^{r\beta_*}}{\Gamma(r\beta_*+1)}.
\]
For fixed $i,j$ the number of chains of length $r$ is at most
$(m-1)^r$. Put $a_R=(m-1)\rho_Rc_{\mathcal B,T}$.
It follows that
\begin{equation}\label{eq:low-S}
 \begin{aligned}
 |S_{ij}(t,\xi)|
 \le\delta_{ij}+\sum_{r=1}^{\infty}
 \frac{a_R^rt^{r\beta_*}}{\Gamma(r\beta_*+1)}=\delta_{ij}+E_{\beta_*,1}(a_Rt^{\beta_*})-1
 \le C_{R,T}.
 \end{aligned}
\end{equation}
This proves absolute and uniform convergence of
\eqref{eq:Sij} on the stated region.

For the source kernels, the same argument gives, for $t>0$,
\[
 \sup_{|\xi|\le R}|S'_{ij}(t,\xi)|
 \le c_{\mathcal B,T}t^{\beta_*-1}
 E_{\beta_*,\beta_*}(a_Rt^{\beta_*}).
\]
Integrating over $(0,T)$, we obtain
\begin{equation}\label{eq:low-Sp}
\int_0^T\sup_{|\xi|\le R}
 |S'_{ij}(s,\xi)|\,\mathrm d s\le c_{\mathcal B,T}\sum_{r=0}^{\infty}
 \frac{a_R^rT^{(r+1)\beta_*}}
 {\Gamma((r+1)\beta_*+1)}=c_{\mathcal B,T}T^{\beta_*}
 E_{\beta_*,\beta_*+1}(a_RT^{\beta_*})
 \le C_{R,T}.
\end{equation}
Thus the series \eqref{eq:Spij} converges absolutely
in $L_1(0,T)$, uniformly for $|\xi|\le R$.

\medskip
\textit{2. The large-frequency region $|\xi|\ge R$.}
By ellipticity and the column-order conditions,
\[
 A_{j,j}(\xi)\asymp|\xi|^{\ell_{jj}},
 \qquad |A_{ij}(\xi)|\le C|\xi|^{\ell_{ij}},
 \qquad |\xi|\ge1.
\]
As in the two-component case, we estimate the coupling terms
relative to the positive diagonal symbols.
For $i\ne j$ put
\[
 q_{ij}(\xi)=
 \frac{\langle \xi\rangle^{\ell_{ii}}|A_{ij}(\xi)|}
 {A_{i,i}(\xi)\langle \xi\rangle^{\ell_{jj}}},
 \qquad q_{ii}(\xi)=0,
 \qquad |\xi|\ge1.
\]
Then
\[
 q_{ij}(\xi)\le C|\xi|^{\ell_{ij}-\ell_{jj}}
 \longrightarrow0\qquad(|\xi|\to\infty).
\]
Since $m$ is fixed, choose $R>1$ so large that
\[
 q_*(\xi):=\max_{1\le i\le m}
 \sum_{j=1}^m q_{ij}(\xi)
 \le\frac12,
 \qquad |\xi|\ge R.
\]
Using $0\le S_j\le1$ and
$\int_0^T K_j(s,\xi)\,\mathrm d s\le1/A_{j,j}(\xi)$,
for each chain in \eqref{eq:Sij} we get
\begin{equation}\label{eq:high-convS}
 (K_{i_0}*\cdots*K_{i_{r-1}}*S_j)(t,\xi)
 \le\prod_{q=0}^{r-1}\frac1{A_{i_q,i_q}(\xi)}.
\end{equation}
Therefore, after multiplying by
$\langle \xi\rangle^{\ell_{ii}-\ell_{jj}}$,
the absolute value of a chain term is bounded by
\[
 \prod_{q=1}^r q_{i_{q-1},i_q}(\xi).
\]
The sum of these products over all chains of length $r$
starting at $i$ is at most $q_*(\xi)^r$.
Thus
\begin{equation}\label{eq:high-S}
 \langle \xi\rangle^{\ell_{ii}-\ell_{jj}}|S_{ij}(t,\xi)|
 \le\sum_{r=0}^{\infty}q_*(\xi)^r\le2,
 \qquad |\xi|\ge R.
\end{equation}
The same bound also proves absolute convergence of the
homogeneous series in the large-frequency region.

For the source kernels, instead of \eqref{eq:high-convS},
we use
\[
 \int_0^T
 (K_{i_0}*\cdots*K_{i_r})(s,\xi)\,\mathrm d s
 \le\prod_{q=0}^r\frac1{A_{i_q,i_q}(\xi)}.
\]
Using \eqref{eq:Spij}, taking absolute values and
summing the chain estimates gives
\begin{equation}\label{eq:high-Sp}
 \langle \xi\rangle^{\ell_{ii}}
 \int_0^T|S'_{ij}(s,\xi)|\,\mathrm d s\le\frac{\langle \xi\rangle^{\ell_{jj}}}{A_{j,j}(\xi)}
 \sum_{r=0}^{\infty}q_*(\xi)^r\le 2\frac{\langle \xi\rangle^{\ell_{jj}}}{A_{j,j}(\xi)}
 \le C_R,
 \qquad |\xi|\ge R.
\end{equation}

For $|\xi|\le R$, all powers $\langle \xi\rangle^\alpha$ are bounded
above and below by positive constants. Consequently,
\eqref{eq:low-S} and \eqref{eq:low-Sp}, combined with
\eqref{eq:high-S} and \eqref{eq:high-Sp}, imply
\eqref{eq:S-bound} and \eqref{eq:Sp-bound} for every
$\xi\in\mathbb R^n$.
\end{proof}

\subsection{Main theorem}\label{sec:main}
\begin{theorem}\label{th:main}
Let $m\ge2$, $\tau>n/2$, and $0<\beta_j\le1$ for
$j=1,\ldots,m$. Assume \eqref{eq:positive} and
\eqref{eq:column}. If
\[
 \varphi_j\in L_2^{\tau+\ell_{jj}}(\mathbb R^n),
 \qquad
 h_j\in C^1([0,T];L_2^\tau(\mathbb R^n)),
 \qquad j=1,\ldots,m,
\]
then the Cauchy problem \eqref{eq:cauchy}--\eqref{eq:initial}
has a unique classical solution
$U\in C([0,T];\mathbf X_\tau)$, given by
\eqref{eq:representation}. Moreover,
\begin{equation}\label{eq:main-regularity}
 \mathbb A(D)U\in C([0,T];\mathbf L_2^\tau),
 \qquad D_t^{\mathcal B}U\in C((0,T];\mathbf L_2^\tau).
\end{equation}
\end{theorem}

\subsection{Proof of Theorem~\ref{th:main}}\label{sec:proof}
\begin{proof}
\textit{1. Existence and continuity.}
By \eqref{eq:S-bound},
\begin{equation}\label{eq:S-op}
 \|S(t,D)\Phi\|_{\mathbf X_\tau}
 \le C_T\|\Phi\|_{\mathbf X_\tau}.
\end{equation}
Consider the source term
\[
 V(t)=\int_0^tS'(t-s,D)H(s)\,\mathrm d s.
\]
Since $H\in C^1([0,T];\mathbf L_2^\tau)$, we have
$H(s)=H(0)+\int_0^sH'(r)\,\mathrm d r$.
For almost every $\xi$ this implies
\begin{equation}\label{eq:source-time-bound}
 |\mathrm F[h_j](s,\xi)|
 \le|\mathrm F[h_j](0,\xi)|
 +\int_0^T|\mathrm F[h_j'](r,\xi)|\,\mathrm d r.
\end{equation}
Taking the Fourier transform of $V$ and using
\eqref{eq:Sp-bound}, we find
\[
\langle \xi\rangle^{\ell_{ii}}|\mathrm F[V_i](t,\xi)|\le C_T\sum_{j=1}^m
 \left(|\mathrm F[h_j](0,\xi)|
 +\int_0^T|\mathrm F[h_j'](r,\xi)|\,\mathrm d r\right).
\]
Multiplying by $\langle \xi\rangle^\tau$, applying Plancherel's theorem
and Minkowski's inequality, we obtain
\[
 \|V(t)\|_{\mathbf X_\tau}
 \le C_T\left(\|H(0)\|_{\mathbf L_2^\tau}
 +\int_0^T\|H'(r)\|_{\mathbf L_2^\tau}\,\mathrm d r\right)
 \le C_T\left(\|H(0)\|_{\mathbf L_2^\tau}
 +T\|H'\|_{C([0,T];\mathbf L_2^\tau)}\right).
\]
Combining this with \eqref{eq:S-op} gives
\[
 \sup_{0\le t\le T}\|U(t)\|_{\mathbf X_\tau}\le C_T\left(\|\Phi\|_{\mathbf X_\tau}
 +\|H(0)\|_{\mathbf L_2^\tau}
 +T\|H'\|_{C([0,T];\mathbf L_2^\tau)}\right).
\]
The same multiplier bounds, together with dominated convergence in
Fourier space, imply
\begin{equation}\label{eq:U-continuity}
 U\in C([0,T];\mathbf X_\tau).
\end{equation}

\medskip
\textit{2. Verification of the equations, initial condition, and regularity.}
The construction of the solution factors and source kernels implies,
for almost every $\xi$ and each $i=1,\ldots,m$,
\begin{equation}\label{eq:full-recurrence}
 \widehat u_i
 =S_i\widehat\varphi_i
 +K_i*\left(\widehat h_i-
 \sum_{\substack{j=1\\j\ne i}}^m A_{ij}(\xi)\widehat u_j\right).
\end{equation}
Indeed, this is precisely the nonhomogeneous version of
\eqref{eq:general-recurrence}, and it also follows term by term from
\eqref{eq:Sij} and \eqref{eq:Spij}.

For the scalar Caputo equation, the relation
$y=S_i y_0+K_i*g$ is equivalent to
\[
 D_t^{\beta_i}y+A_{ii}(\xi)y=g,
 \qquad y(0)=y_0,
\]
as follows from \eqref{eq:scalar-L}. Applying this identity to
\eqref{eq:full-recurrence} gives
\begin{equation}\label{eq:fourier-equation-verified}
 D_t^{\mathcal B}\widehat U(t,\xi)
 +\mathbb A(\xi)\widehat U(t,\xi)
 =\widehat H(t,\xi),
 \qquad t>0,
\end{equation}
for almost every $\xi$.

Next, by \eqref{eq:column}, $\ell_{ij}\le\ell_{jj}$ for all $i,j$.
Hence
\[
 \|A_{ij}(D)u_j(t)\|_{L_2^\tau}
 \le C\|u_j(t)\|_{L_2^{\tau+\ell_{ij}}}
 \le C\|u_j(t)\|_{L_2^{\tau+\ell_{jj}}}.
\]
Together with \eqref{eq:U-continuity}, this yields
\begin{equation}\label{eq:A-continuity}
 \mathbb A(D)U\in C([0,T];\mathbf L_2^\tau).
\end{equation}
Therefore \eqref{eq:fourier-equation-verified} can be inverted in
$\mathbf L_2^\tau$, and the equation \eqref{eq:cauchy} holds for
$t>0$. Moreover,
\[
 D_t^{\mathcal B}U=H-\mathbb A(D)U
 \in C((0,T];\mathbf L_2^\tau),
\]
which proves \eqref{eq:main-regularity}.

The bounds \eqref{eq:S-bound}--\eqref{eq:Sp-bound} and
\eqref{eq:source-time-bound} provide integrable Fourier majorants.
Thus dominated convergence yield
\[
U(t)\longrightarrow\Phi\quad\text{in }\mathbf X_\tau
 \qquad(t\rightarrow0).
\]
Hence \eqref{eq:initial} holds.

Finally, since $\tau>n/2$, \eqref{eq:embedding},
\eqref{eq:U-continuity}, and \eqref{eq:A-continuity} imply the spatial
continuity required in \eqref{eq:classical}; the same conclusion for
$D_t^{\mathcal B}U$ follows from the equation. Thus the function
\eqref{eq:representation} is a classical solution.

\medskip
\textit{3. Uniqueness.}
The uniqueness is proved in the same way as
in the two-component case. Indeed, the
difference of two solutions with the same
data satisfies a homogeneous system with
zero initial conditions. Applying the
Fourier transform, we obtain, for almost
every fixed $\xi$, a linear system of
fractional ordinary differential equations
with constant coefficients and zero initial
data. By the standard uniqueness result
for linear multi-order fractional systems
(see \cite{UmarovNew}), the Fourier
transform of the difference vanishes.
Consequently, the two solutions coincide.
This completes the proof.
\end{proof}

\section{Examples}\label{sec:examples}

Let $\tau>n/2$, $h_j\in C^1([0,T];L_2^\tau(\mathbb R^n))$, and
$u_j(0,x)=\varphi_j(x)$. All convolutions are taken in $t$ at fixed
$\xi$ and act entrywise on matrices, with $f^{*0}*g=g$.

\subsection*{Example 1. A two-component diffusion system with rotational coupling}

Let $\kappa,\omega>0$ and $0<\beta_1,\beta_2\le1$. Consider
\[
\begin{cases}
D_t^{\beta_1}u_1-\kappa\Delta u_1+\omega u_2=h_1,\\
D_t^{\beta_2}u_2-\omega u_1-\kappa\Delta u_2=h_2 .
\end{cases}
\]
For $\varphi_1,\varphi_2\in L_2^{\tau+2}(\mathbb R^n)$,
Theorem~\ref{th:two} gives a unique classical solution in
$C([0,T];(L_2^{\tau+2})^2)$,
\[
 u_i(t,x)=\sum_{j=1}^{2}\left[S_{ij}(t,D)\varphi_j(x)
 +\int_0^t S'_{ij}(t-s,D)h_j(s,x)\,\mathrm d s\right],\quad i=1,2,
\]
where
\[
\begin{aligned}
(S_{ij})&=\sum_{k=0}^{\infty}(-\omega^2)^k(K_1*K_2)^{*k}*
\begin{pmatrix}
 S_1 & -\omega(K_1*S_2)\\
 \omega(K_2*S_1)&S_2
\end{pmatrix},\\
(S'_{ij})&=\sum_{k=0}^{\infty}(-\omega^2)^k(K_1*K_2)^{*k}*
\begin{pmatrix}
 K_1&-\omega(K_1*K_2)\\
 \omega(K_2*K_1)&K_2
\end{pmatrix},
\end{aligned}
\]
with
\[
 S_j=E_{\beta_j}(-\kappa|\xi|^2t^{\beta_j}),\qquad
 K_j=t^{\beta_j-1}E_{\beta_j,\beta_j}(-\kappa|\xi|^2t^{\beta_j}),
 \quad j=1,2.
\]

\subsection*{Example 2. A fully coupled system with three independent time orders}

Let $a,b,c>0$ and $0<\beta_1,\beta_2,\beta_3\le1$. Consider
\[
\begin{cases}
D_t^{\beta_1}u_1-\Delta u_1+a\partial_{x_1}u_2+b\Delta u_3=h_1,\\
D_t^{\beta_2}u_2+a\partial_{x_1}u_1-\Delta u_2-b\Delta u_3=h_2,\\
D_t^{\beta_3}u_3-c(u_1+u_2)+\Delta^2u_3=h_3 .
\end{cases}
\]
Its Fourier symbol is
\[
\mathcal A(\xi)=
\begin{pmatrix}
 |\xi|^2 & \mathrm i a\xi_1 & -b|\xi|^2\\
 \mathrm i a\xi_1 & |\xi|^2 & b|\xi|^2\\
 -c & -c & |\xi|^4
\end{pmatrix}.
\]
The diagonal orders are $2,2,4$, and the largest off-diagonal orders
in the corresponding columns are $1,1,2$. Thus for
$\Phi\in\mathbf X_\tau=(L_2^{\tau+2})^2\times L_2^{\tau+4}$,
Theorem~\ref{th:main} gives a unique classical solution in
$C([0,T];\mathbf X_\tau)$,
\[
 u_i(t,x)=\sum_{j=1}^{3}\left[S_{ij}(t,D)\varphi_j(x)
 +\int_0^t S'_{ij}(t-s,D)h_j(s,x)\,\mathrm d s\right],\quad i=1,2,3,
\]
where $S_{ij}=I_t^{1-\beta_j}S'_{ij}$, $I_t^0=I$, and
\[
\begin{aligned}
(S'_{ij})={}&\sum_{k=0}^{\infty}
\Bigl[-a^2\xi_1^2(K_1*K_2)
 +bc|\xi|^2(K_1*K_3-K_2*K_3)\Bigr]^{*k}\\
&*\left[\begin{pmatrix}
 K_1 & -\mathrm i a\xi_1(K_1*K_2) & b|\xi|^2(K_1*K_3)\\
 -\mathrm i a\xi_1(K_1*K_2) & K_2 & -b|\xi|^2(K_2*K_3)\\
 c(K_1*K_3)&c(K_2*K_3)&K_3
\end{pmatrix}\right.\\
&\hspace{1em}\left.{}+(K_1*K_2*K_3)
\begin{pmatrix}
 bc|\xi|^2&bc|\xi|^2&\mathrm i ab\xi_1|\xi|^2\\
 -bc|\xi|^2&-bc|\xi|^2&-\mathrm i ab\xi_1|\xi|^2\\
 -\mathrm i ac\xi_1&-\mathrm i ac\xi_1&a^2\xi_1^2
\end{pmatrix}\right].
\end{aligned}
\]
Here
\[
\begin{aligned}
 K_j&=t^{\beta_j-1}E_{\beta_j,\beta_j}(-|\xi|^2t^{\beta_j}),
 &&j=1,2,\\
 K_3&=t^{\beta_3-1}E_{\beta_3,\beta_3}(-|\xi|^4t^{\beta_3}).
\end{aligned}
\]

\section*{Conclusion}

In this article, we established sufficient conditions for the classical
solvability of the Cauchy problem for linear multi-order Caputo systems
with full matrices of constant-coefficient spatial differential operators.
The construction does not require rational relations among the fractional
orders or a triangular spatial operator matrix. By combining Fourier--Laplace
representations with series of time convolutions of scalar Mittag--Leffler
kernels, we obtained solution and source operators that account for repeated
interactions between components.

\section*{Author contributions}

All authors contributed equally to this work.

\section*{Conflict of interest}

The authors declare no conflict of interest.


\begin{center}
References
\end{center}
\begin{enumerate}
\setlength{\itemsep}{3pt}

\bibitem{Umarovappear}
S.~Umarov,
Representations of solutions of systems of time-fractional
pseudo-differential equations,
\emph{Fractional Calculus and Applied Analysis}
\textbf{27} (2024), no.~2, 616--651.
{doi:10.1007/s13540-024-00241-z}.

\bibitem{UmarovNew}
S.~Umarov,
Representations of solutions of time-fractional multi-order systems
of differential-operator equations,
\emph{Fractal and Fractional} \textbf{8} (2024), no.~5, Article~254.
{doi:10.3390/fractalfract8050254}.

\bibitem{AO}
R.~Ashurov and O.~Mukhiddinova,
On systems of fractional nonlinear partial differential equations,
\emph{Fractional Calculus and Applied Analysis}
\textbf{28} (2025), no.~5, 2167--2197.
{doi:10.1007/s13540-025-00431-3}.

\bibitem{AS}
R.~Ashurov and I.~Sulaymonov,
Fractional parabolic systems of vector order,
\emph{Journal of Mathematical Sciences}
\textbf{284} (2024), no.~2, 179--195.
{doi:10.1007/s10958-024-07342-3}.

\bibitem{SulaymonovRL}
I.~Sulaymonov,
Fractional parabolic systems of vector order with Riemann--Liouville
fractional derivative,
in R.~Ashurov and E.~Karimov (eds.),
\emph{Modern Methods in Mathematical Physics and their Applications},
Trends in Mathematics, vol.~16, Birkh\"auser, Cham, 2026, pp.~277--284.
{doi:10.1007/978-3-032-13872-9\_32}.



\bibitem{M}
F.~Mainardi,
On some properties of the Mittag--Leffler function
$E_\alpha(-t^\alpha)$, completely monotone for $t>0$ with $0<\alpha<1$,
\emph{Discrete and Continuous Dynamical Systems, Series B}
\textbf{19} (2014), no.~7, 2267--2278.
{doi:10.3934/dcdsb.2014.19.2267}.

\bibitem{Gor}
R.~Gorenflo, A.~A.~Kilbas, F.~Mainardi, and S.~V.~Rogosin,
\emph{Mittag-Leffler Functions, Related Topics and Applications},
Springer, Berlin, Heidelberg, 2014.
{doi:10.1007/978-3-662-43930-2}.

\bibitem{Umarov1}
S.~Umarov,
On fractional Duhamel's principle and its applications,
\emph{Journal of Differential Equations} \textbf{252} (2012),
no.~10, 5217--5234.
{doi:10.1016/j.jde.2012.01.029}.

\bibitem{Umarov2}
S.~Umarov,
\emph{Introduction to Fractional and Pseudo-Differential Equations
with Singular Symbols}, Springer, Cham, 2015.
{doi:10.1007/978-3-319-20771-1}.

\end{enumerate}
\end{document}